\documentclass[11 pt, reqno]{article}

\usepackage{amsmath,amsfonts,amssymb,amsthm}
\usepackage[colorlinks=true,hyperindex=true]{hyperref}
\usepackage{cancel,bbm}
\usepackage{mathabx} 
\usepackage{comment}
\usepackage{enumitem}
\usepackage{cite}

\usepackage{chngcntr}
\counterwithin*{equation}{section}

\usepackage[left=3cm, right=3cm, top=3 cm, bottom= 3 cm]{geometry}
\usepackage{color}
\usepackage{fancyhdr}
\usepackage{latexsym}

\usepackage{bm}

\newtheorem{ccounter}{ccounter}[section]
\newtheorem{thm}[ccounter]{Theorem}
\newtheorem{lem}[ccounter]{Lemma}
\newtheorem{cor}[ccounter]{Corollary}
\newtheorem{defn}[ccounter]{Definition}
\newtheorem{prop}[ccounter]{Proposition}
\newtheorem{ass}[ccounter]{Assumption}
\newtheorem{ex}[ccounter]{Example}

\NewCommandCopy{\dotlessi}{\i}

\def\bet{\begin{thm}}
\def\eet{\end{thm}}
\def\bel{\begin{lem}}
\def\eel{\end{lem}}
\def\bas{\begin{ass}}
\def\eas{\end{ass}}
\def\bec{\begin{cor}}
\def\eec{\end{cor}}
\def\bed{\begin{defn}}
\def\eed{\end{defn}}
\def\bep{\begin{prop}}
\def\eep{\end{prop}}
\def\beq{\begin{equation}}
\def\eeq{\end{equation}}
\def\proof{\noindent {\bf Proof.}\ \ }
\def\bea{\begin{equation*}}
\def\eea{\end{equation*}}

\def\bex{\begin{ex}}
\def\eex{\end{ex}}

\def\rr{\mathbb{R}}

\def\cc{\mathbb{C}}
\def\1{\boldsymbol{1}}
\def\Im{\mathrm{Im}}

\def\e{\mathrm{e}}
\def\i{\mathrm{i}}
\def\del{\partial}
\def\d{\mathrm{d}}
\def\eps{\varepsilon}
\renewcommand\leq\varleq
\renewcommand\geq\vargeq
\def\ee{\mathrm{E}}

\def\F{\mathcal{F}}
\def\O{\mathcal{O}}

\def\ee{\mathbb{E}}

\def\pp{\mathbb{P}}

\def\msc{m_{\mathrm{sc}}}
\def\rhosc{\rho_{\mathrm{sc}}}

\def\mfa{\mathfrak{m}}
\def\A{\mathcal{A}}

\def\bx{\mathbf{x}}
\def\by{\mathbf{y}}

\def\mfa{\mathfrak{a}}

\def\mfb{\mathfrak{b}}

\def\haty{\hat{y}}
\def\B{\mathcal{B}}

\def\bfv{\boldsymbol{v}}
\def\bx{\boldsymbol{x}}
\def\bz{\boldsymbol{z}}
\def\nn{\mathbb{N}}
\def\N{\mathcal{N}}
\def\nsc{n_{\mathrm{sc}}}
\def\tilgam{\tilde{\gamma}}
\def\hatgam{\hat{\gamma}}

\def\haty{\hat{y}}
\def\hfb{\hat{\mathfrak{b}}}
\def\tily{\tilde{y}}

\begin{document}

\begin{table}
\centering

\begin{tabular}{c}

\multicolumn{1}{c}{\parbox{12cm}{\begin{center}\Large{\bf Multilevel limits of spiked random matrix corners processes}\end{center}}}\\
\\
\end{tabular}
\begin{tabular}{ c c c  }
Benjamin Landon
& \phantom{blah} & 
Tara Stojimirovic
 \\
 & & \\  
 \small{University of Toronto} & & \small{University of Toronto} \\
 \small{Department of Mathematics} & & \small{Department of Mathematics} \\
 \small{\texttt{blandon@math.toronto.edu}} & & \small{\texttt{tara.stojimirovic@mail.utoronto.ca}} \\
  & & \\
\end{tabular}
\\
\begin{tabular}{c}
\multicolumn{1}{c}{\today}\\
\\
\end{tabular}

\begin{tabular}{p{15 cm}}
\small{{\bf Abstract:} We consider the largest eigenvalues of submatrices of the classical Gaussian random matrix ensembles with a finite rank spike in the critical BBP regime. We show that they converge to a multilevel interlacing particle system whose lowest level is the Airy$_\beta$ point process. Our construction and scaling limits extend in a natural way to non-classical $\beta >0$. 
}
\end{tabular}
\end{table}

\section{Introduction}

Let $G$ be an $(N+K) \times (N+K)$ matrix from the Gaussian $\beta$-ensemble for the classical values $\beta \in \{1, 2, 4\}$, scaled so that $\ee[ |G_{ij} |^2 ] =1$ for $i \neq j$. Let $A$ be an $(N+K) \times (N+K)$  matrix that is $0$ everywhere except for the bottom right $K$ diagonal entries where
\beq
A_{(N+i), (N+i)} = N^{1/2} + \alpha_i N^{1/6}, \qquad 1 \leq i \leq K
\eeq
for some  $\alpha_i \in \rr$. We will think of $K$ as fixed and consider the limit $N \to \infty$. When $\beta=2$ the scaling limit of the largest eigenvalues of $G+A$ was first found in the celebrated work of Baik, Ben Arous and P{\'e}ch{\'e} \cite{baik2005phase} (at least for the sample covariance matrix version of the problem; the additive Hermitian model was handled in \cite{peche2006largest} via similar techniques). In the case $K=1$, the  limiting distribution of the largest eigenvalue interpolates between the Tracy-Widom and Gaussian distributions as $\alpha $ ranges from $- \infty$ to $+\infty$,  a phenomenon  now known as the BBP transition. We emphasize that our choice of scaling places our model in the critical regime where the asymptotic fluctuations are neither Tracy-Widom nor Gaussian. 

The $\beta=2$ case relied on the fact that the correlation functions of the eigenvalues of $G+A$ are determinantal. In particular, these methods did not extend to the cases $\beta =1, 4$. Bloemendal and Vir{\'a}g gave a new approach based on tri-diagonal random matrix models \cite{alex2013limits,bloemendal2016limits}, associating the scaling limits of the largest eigenvalue of $G+A$ with those of certain random differential operators.

In fact, when $\beta=2$ (and for  $A$ of arbitrary rank)  the eigenvalues of $G+A$ and those of its leading principal  submatrices\footnote{A principal submatrix is a square submatrix retaining the same set of row and column indices. A leading principal submatrix is formed by the top left corner of the original matrix.} together form a determinantal point process  on $\mathbb{N} \times \rr$ \cite{adler2013random}. Using the determinantal structure,   the largest eigenvalues of the top $K+1$ leading principal submatrices (when $A$ has rank $K$) have been shown to converge to a multilevel interlacing particle system whose bottom level is the Airy$_2$ point process \cite{bao2022eigenvector}.

In this paper we will extend this result to $\beta=1, 4$ and also in a natural way to general $\beta >0$, where there is no known random matrix model giving the interlacing system in finite dimensions. 

Let us explain our result in the simplest case when $K=1$ and $\beta \in \{1, 2, 4\}$. Then $A$ has only a single non-zero entry in the bottom right corner, which we denote $N^{1/2} + \alpha N^{1/6}$. Let $ \lambda_1 > \lambda_2 > \dots $ be the eigenvalues of the top left $N \times N$ submatrix of $G$ with $\ell^2$-normalized eigenvectors $\{u_i\}_{i=1}^N$.  Denote the right-most column of $G$ as $( \bfv^T , g )^T$, with $\bfv$ an $N$-dimensional vector with values in $\rr, \cc$ or $\mathbb{H}$ depending on the value of $\beta$. Let $\mu_1 > \mu_2 > \dots$ denote the eigenvalues of $G+A$ which satisfy the interlacing condition $\mu_1 > \lambda_1 > \mu_2 > \lambda_2 > \dots > \lambda_N > \mu_{N+1}$, almost surely. Then by elementary linear algebraic manipulations the $\{ \mu_i\}_{i=1}^{N+1}$ are the $N+1$ unique solutions to the equation (see, e.g., \cite[(2.3)]{huang2019eigenvalues})
\beq \label{eqn:intro-3}
 g+ N^{1/2} + \alpha N^{1/6}= \mu + \sum_{i=1}^N \frac{ | u_i^* \bfv  |^2}{ \lambda_i - \mu } =: \mu + \sum_{i=1}^N \frac{ z_i}{ \lambda_i - \mu }
\eeq
Due to the rotational invariance of $\bfv$ and its independence from the eigenbasis $\{ u_i \}_{i=1}^N$ we have that $z_i := | u_i^* \bfv   |^2$ are independent random variables with Gamma distribution having density proportional to $x^{\beta/2-1} \e^{- \beta x /2} \1_{ \{ x > 0 \} }$. Note that $\ee[z_i] =1$. 

If we define 
\beq \label{eqn:GBE-rescaled}
x_i := N^{1/6} \left( 2 N^{1/2} - \lambda_i \right)
\eeq
then the $x_i$ converge in the finite dimensional distribution sense to a point process $ \mfa_1 < \mfa_2 < \dots$ known as the Airy$_\beta$ point process \cite{ramirez2011beta,tracy1994level,tracy1996orthogonal}. In analogy with the BBP limits, we introduce the re-scaled quantities $y_i := N^{1/6} ( 2 N^{1/2} - \mu_i )$ in order to investigate the scaling limits of the $\mu_i$. By manipulating \eqref{eqn:intro-3} we see that the $\{ y_i \}_{i=1}^{N+1}$ solve the equation 
\beq \label{eqn:intro-1}
\alpha = \sum_{i=1}^N \frac{z_i}{y -x_i} + N^{1/3} - \left( y N^{-1/3} + g N^{-1/6} \right). 
\eeq
If we recognize the diverging quantity $N^{1/3}$ as the contribution from the limiting semicircle distribution, then it is natural to expect that the $y_i$ converge to a point process $\mfb_i$ that is defined as the unique solutions (in the variable $y$) to the equation,
\beq \label{eqn:intro-2}
\alpha = \lim_{n \to \infty} \left( \sum_{i=1}^n \frac{z_i}{ y - \mfa_i } + \int_0^{\left(\frac{3\pi n}{2} \right)^{2/3}} \frac{ x^{-1/2}}{\pi} \d x \right) ,
\eeq
satisfying $\mfb_1 < \mfa_1 < \mfb_2 < \mfa_2 < \dots $. We ignore the question about existence of the limit on the RHS of \eqref{eqn:intro-2} for the present discussion. We remark that the deterministic term on the RHS  is the natural renormalization to the diverging sum (in our scaling we expect $\mfa_i \approx i^{2/3}$ and $\sum_{i \geq 1} i^{-2/3}$ is divergent). 

The joint convergence of the $\{ x_i\}_i$ and $\{ y_i\}_i$ (in the finite dimensional sense) to the point processes $\{\mfa_i\}_i$ and $\{ \mfb_i\}_i$ described above is part of what we prove. Our main result extends this in two natural directions. First, it is clear that this discussion extends in a natural way to multiple levels. That is, the eigenvalues of the $(N+j) \times (N +j)$ leading principal submatrix of $G+A$ depend in a similar  Markovian way on the eigenvalues of the $(N+j-1) \times (N+j-1)$ leading principal  submatrix as to what has been described above, which is the case $j=1$. One need only introduce a new set of independent $z_i$ in \eqref{eqn:intro-3}. 
Indeed, the new set of $z_i$ come from the topmost $N+j-1$ entries of the $(N+j)$th column of $G$ and so always have the same distribution no matter the choice of $A$, and are independent of the $(N+j-1) \times (N+j-1)$  leading principal submatrix of $G+A$. Using this observation, we find a  joint limit of the largest eigenvalues of the top $K+1$ principal submatrices to a sequence of interlacing point processes where the $j$th level $\{\mfb_i\}_i$ is constructed from the $(j-1)$st level $\{\mfa_i\}_i$ as in \eqref{eqn:intro-2} and the first level is the Airy$_\beta$ point process.

Secondly, the above discussion extends immediately to non-classical $\beta >0$ in a natural way. That is, if $\lambda_1 > \dots > \lambda_N$ are the eigenvalues of the $N$-dimensional Gaussian $\beta$-ensemble, then given a doubly indexed sequence of independent $z_i^{(j)}$ having density proportional to $x^{\beta/2-1} \e^{ - \beta x/2}\1_{ \{ x > 0 \} }$ and centered Gaussians $g^{(j)}$ with variance $2/\beta$, one can successively define a sequence of point processes via \eqref{eqn:intro-3} (equivalently \eqref{eqn:intro-1} after rescaling). The scaling limit to a multilevel system of interlacing particles is described in an identical fashion as in the case of classical $\beta$.

In fact when  $A =0$, the prelimiting process is known as the Hermite $\beta$-corners process and was introduced by Gorin and Shkolnikov with an explicit density in \cite{gorin2015multilevel}. To our knowledge, the connection between \cite{gorin2015multilevel} and the sequence of solutions to \eqref{eqn:intro-3} was first pointed out in \cite{najnudel2021bead}; this relation follows from \cite[Proposition 4.3.2]{forrester2010log} which computes the density of zeros of the RHS of \eqref{eqn:intro-1} conditional on the $x_i$. In fact, a simple change of variable in \cite[Proposition 4.3.2]{forrester2010log} reveals that for general diagonal $A$, the eigenvalues of the submatrices of $G+A$ have an explicit joint density. This process has recently been studied in \cite{petrov2026perturbedbetacornersprocess}, which appeared on the arXiv while we were preparing this paper. Our work can be interpreted as finding a scaling limit of the process considered in \cite{petrov2026perturbedbetacornersprocess} for finite rank $A$ in the critical BBP regime.

Finally let us comment on the choice of scaling of the perturbation: that is, why we choose $\alpha_i \in \rr$ fixed in the definition of $A$. If, say, we have that $\alpha_i \to - \infty$ as $N \to \infty$ for some $i$, then we expect that the difference between  the $j$th particles in the consecutive levels associated to that $\alpha_i$ tends to $0$ as $N\to \infty$. This can be read off at a formal level from \eqref{eqn:intro-2}. In this sense, this choice is degenerate. One could scale the difference between the $j$th particles in consecutive levels, finding a distribution coming from the $z_i$. This has been done in the $A=0$ case in \cite{huang2019eigenvalues,gorin2014interactingparticlesystemsedge}, and we expect the more general case to be similar. On the other hand if $\alpha_i \to + \infty$ then one has a similar degeneration. The $j$th particle in the lower level sticks to the $(j+1)$st particle in the higher level. The smallest particle in the higher level would converge to $-\infty$ in our scaling, and we expect that, under an appropriate rescaling, it would have Gaussian fluctuations independent of the rest of the particles, as in the supercritical BBP regime. We therefore require $\alpha_i$ to remain finite in order to have the nontrivial multilevel structure given by \eqref{eqn:intro-2}.

\subsection{Statement of main result} \label{sec:mr}

\subsubsection{Definition of processes} \label{sec:def-processes}

We now formalize the discussion in the previous section. We first define a map $\Phi_{N, i}$ that outputs the solutions to the equation \eqref{eqn:intro-1} given the other data. That is, we let $N \geq 1$, $i \geq 0$ be integers. Then $\Phi_{N, i}$ takes as input vectors $\bx \in \rr^{N+i}$, $\bz \in (0, \infty)^{N+i}$, and scalars $g, \alpha \in \rr$. We also assume that the entries of $\bx = (x_1, \dots, x_{N+i})$ are in increasing order, $x_1 < x_2 < x_3 < \dots $. Given this data, the map $\Phi_{N, i} : \rr^{N+i} \times (0, \infty)^{N+i} \times \rr \times \rr \to \rr^{N+i+1}$ outputs the $N+i+1$ solutions $\by = (y_1, \dots y_{N+i+1} ) \in \rr^{N+i+1}$ to the equation
\beq \label{eqn:yhat-def}
\alpha = \sum_{j=1}^{N+i} \frac{z_j}{y-x_j} + N^{1/3} - (y N^{-1/3} + g N^{-1/6} )
\eeq
labelled in increasing order $y_1 < y_2 < \dots < y_{N+i+1}$. By the discussion in the previous section, if:
\begin{itemize}
\item $x_1 < x_2 <  \dots < x_N$ are the re-scaled eigenvalues of the Gaussian $\beta$-ensemble
\end{itemize}
and for every $i \geq 1$,
\begin{itemize}
\item  $\bz^{(i)} = ( z_1^{(i)} , z_2^{(i)} , \dots )$ has entries which are i.i.d. random variables with density 
\beq
\pp\left[ z  \in \d x \right] =  \frac{1}{C_\beta} x^{\beta/2-1} \e^{ - \beta  x /2 } \1_{ \{ x > 0 \} } \d x
\eeq
where $C_\beta >0$ is a normalization constant
\item $g_i$ are i.i.d. centered Gaussians with variance $2 / \beta$
\item $\alpha_i \in \rr$
\item All of the random variables $z_j^{(i)}$ and $g_i$   are mutually independent for all $i, j$ and also of the family $\{ x_k \}_{k=1}^N$
\end{itemize}
then the joint distribution of the sequence defined inductively by 
$$\bx^{(i)} = \Phi_{N, i-1} ( \bx^{(i-1)} , (z_j^{(i)} )_{j=1}^{N+i-1} , g_i, \alpha_i)$$ is that of the re-scaled  eigenvalues of $G+A$ and its leading principal submatrices in the classical cases $\beta=1, 2, 4$. Here we denoted $\bx^{(0)} = (x_1 ,  \dots , x_N )$. However, we consider the above construction for any $\beta >0$ as it makes no reference to whether or not $\beta \in \{1, 2, 4\}$.

We now define the scaling limit. Given a choice of countably many distinct $\mfa = ( \mfa_1, \mfa_2 , \dots ) \in \rr^{\nn}$  labelled in increasing order and $\bz \in (0, \infty)^\nn$ let us define the function,
\beq \label{eqn:F-def}
F (E ; \mfa, \bz) := \lim_{n \to \infty} \left( \sum_{i=1}^n \frac{z_i}{ E - \mfa_i } + \int_0^{\left( \frac{3 n \pi}{2} \right)^{2/3}} \frac{ x^{-1/2}}{ \pi} \d x \right).
\eeq
We say $(\mfa, \bz)$ is \emph{admissible} if the limit exists for every $E \in \rr \backslash \{ \mfa \}$ and defines a strictly decreasing continuous function on the connected components of $\rr \backslash \{ \mfa \}$ s.t. for all $i \geq 1$,
\beq
\lim_{ E \to \mfa_i^-} F(E ; \mfa, \bz ) = - \infty, \qquad \lim_{E \to \mfa_{i-1}^+} F(E ; \mfa, \bz ) = + \infty
\eeq
where we defined $\mfa_0 := - \infty$. If the pair $(\mfa, \bz)$ is admissible, then it is clear that we can define a process $\mfb \in \rr^\nn$, labelled in increasing order as the unique solutions to
\beq \label{eqn:mfb-def}
\alpha = \lim_{n \to \infty} \left( \sum_{i=1}^n \frac{z_i}{ E - \mfa_i} + \int_0^{  \left( \frac{ 3 n \pi}{2} \right)^{2/3} }\frac{1}{ \pi x^{1/2}} \d x \right) = F(E ; \mfa , z)
\eeq
that interlace $\mfa$ in the sense that $\mfb_1 < \mfa_1 < \mfb_2 < \mfa_2 < \dots $.  We define $\mfb = \Phi ( \mfa, \bz, \alpha )$. Here $\Phi : \rr^{\nn} \times (0, \infty)^{\nn} \times \rr \to \rr^{\nn}$ is defined on triples $(\mfa, \bz , \alpha ) \in \rr^{\nn} \times ( 0 , \infty )^{\nn} \times \rr$ s.t. $(\mfa, \bz)$ is admissible.
 The following states that our candidate scaling limit is well-defined and is proven in Section \ref{sec:existence}.

\bep \label{prop:intro-well-defined}
Let $\mfa$ be the Airy $\beta$-ensemble, let $\alpha_i \in \rr$ and let $\bz^{(i)}$ be independent vectors of i.i.d.  random variables with density proportional to $x^{\beta/2-1} \e^{ - \beta x /2} \1_{ \{ x > 0 \} }$, independent of $\mfa$. Then $(\mfa, \bz^{(1)})$ is almost surely admissible and the sequence defined inductively by $\mfa^{(i)} = \Phi ( \mfa^{(i-1)} , \bz^{(i)} , \alpha_i)$, where we denote $\mfa^{(0)} = \mfa$, has the property that $( \mfa^{(i)} , \bz^{(i+1)})$ is almost surely admissible for all $i$.  
\eep

\subsubsection{Statement of main result}

Our main result is the following:

\bet \label{thm:main}  For any $K \geq 1$, the sequence $\{ \bx^{(i)}\}_{i=0}^{K}$ defined in Section \ref{sec:def-processes} converges to $\{\mfa^{(i)}\}_{i=0}^{K}$ in the sense of finite dimensional distributions. 
\eet

This theorem is proved in Section \ref{sec:convergence}. Let us comment now briefly on our methods, which are quite natural. We focus on the $K=1$ case, the general case being similar. The main challenge is that we only know the finite dimensional convergence of the $x_i$ to $\mfa_i$, and \eqref{eqn:yhat-def} and \eqref{eqn:mfb-def} involve a growing number of particles.

We therefore consider truncated analogs of \eqref{eqn:yhat-def} and \eqref{eqn:mfb-def}. That is, instead of considering \eqref{eqn:yhat-def} with $N$ particles and \eqref{eqn:mfb-def} with infinitely many particles, we consider both with only a fixed number of particles $m$. We must show that for a large fixed $m$, the solutions of the truncated systems approximate the true solutions with high probability, uniformly as $N \to \infty$. This requires a careful handling of the tail of the sum in \eqref{eqn:yhat-def}. In particular, since $m$ must be fixed, this puts the estimates out of reach of the well-known rigidity results in random matrix theory \cite{erdHos2017dynamical}. Instead we prove tightness results for the largest particles of the Gaussian $\beta$-ensemble using the work of \cite{bourgade2022optimal}. While  the form of tightness we obtain is quite strong, weaker results would suffice. For example, the alternative estimates of \cite{landon2022fluctuations} derived for the classical values of $\beta$ (relying in particular on Gustavsson's classic paper on the GUE \cite{gustavsson2005gaussian}) would suffice.

\subsection{Relation to other work}

Some of our methods are similar to those of Bykhovskaya, Gorin and Sodin \cite{bykhovskaya2026weakweakfactorsuniform}. That paper considers the $K=1$ case and shows that the largest eigenvalues of $G+A$ scale to the solutions of \eqref{eqn:intro-2}, giving an alternative characterization of the distributions found in \cite{alex2013limits}. In the case that $\mfa_i$ is the Airy point process they dub the RHS of \eqref{eqn:intro-2} the ``Airy Green function,'' as it coincides with the scaling limit of a resolvent entry of the classical G$\beta$E. The main focus of that paper is on the statistical applications of the BBP transition. In particular, the work \cite{bykhovskaya2026weakweakfactorsuniform}  did not investigate the possibilities of constructing a multilevel scaling limit as we do here. Furthermore, their work finds the limit by a  different method: it uses rank one perturbation theory to write the eigenvalues of $G+A$ in terms of $G$ instead of the leading principal submatrix. Due to the aforementioned eigenvalue sticking phenomena (that is, the behavior when $\alpha_i \to - \infty$ discussed just before Section \ref{sec:mr}), this leads to the same limits. At a more technical level,  their proof of convergence is different from ours, as they rely more heavily on meromorphicity properties of random matrix Stieltjes transforms; ours is more rudimentary, and so may be of relevance  in other settings. Additionally, they only consider the classical $\beta$, as these are the most interesting cases for the statistical applications of their paper. Extending their work to general $\beta >0$ requires adequate tightness estimates for the Gaussian $\beta$-ensemble. We derive these from the work \cite{bourgade2022optimal}. In fact, our work verifies two conjectures of \cite[Remarks 4.5 and 4.10]{bykhovskaya2026weakweakfactorsuniform} in the Gaussian setting considered here: the extension to general $\beta >0$ and the characterization of the higher rank spiked distributions as solutions to equations involving \eqref{eqn:mfb-def}.

The work \cite{bao2026eigenvectordistributionrandommatrices} also considers the Airy Green function of \cite{bykhovskaya2026weakweakfactorsuniform}. They use it to find the distribution of eigenvector components of spiked random matrices in the critical BBP regime, under general finite rank perturbations (the case $\beta=2$ being handled in \cite{bao2022eigenvector} by partially relying on determinantal methods). Their limit involves a derivative of the Airy Green function. They also handle the rank 1 critically spiked case for general $\beta >0$. Part of their work extends \cite{bykhovskaya2026weakweakfactorsuniform} to general finite rank $K$, describing the limit of the eigenvalues of $G+A$ in a similar fashion to \cite{bykhovskaya2026weakweakfactorsuniform} (using similar notions of convergence to \cite{bykhovskaya2026weakweakfactorsuniform} which, as we stated, are different from ours). As their work has a different focus -- the eigenvector components -- they did not consider the multilevel joint limits as we do here. Moreover, their work uses convergence of the spiked eigenvalues of Bloemendal and Vir{\'a}g \cite{alex2013limits,bloemendal2016limits} as input, which we do not rely on.

Similar scaling limits have been derived in the context of the $2$-spin spherical SK model \cite{landon2022fluctuations,landon2020fluctuations} and some of the methods here are similar to the methods there. 

As we mentioned above, while we were preparing this paper, the aforementioned \cite{petrov2026perturbedbetacornersprocess} appeared. It studies the processes we consider but for fixed $N$, and focuses on various asymptotic regimes as $\beta \to \infty$. 

The interlacing system we derive could be viewed as an edge analog of \cite{najnudel2021bead} that derives the bulk limit of the Hermite $\beta$-corners process (i.e., with $A=0$) as a bead process involving the $\mathrm{Sine}_\beta$ process. In particular, they derive the version of \eqref{eqn:mfb-def} in the bulk (i.e., the $j$th level of the system is derived from the $(j-1)$st level via  a similar Markovian structure involving a Stieltjes transform). Note that without a spike, the limiting system \eqref{eqn:mfb-def} is degenerate in the sense described above at the edge of the spectrum; the $j$th level of particles just equals the $(j-1)$st level. This is in contrast to the bulk setting which does not require a spike for non-degeneracy.

\vspace{5 pt}

\noindent{\bf Acknowledgements}. The work of the authors is supported by NSERC. BL thanks Vadim Gorin for helpful discussions.

\vspace{5 pt}

\noindent{\bf AI Disclosure.} The mathematical arguments and proofs were conceived and written by the authors without AI assistance. After uploading the first two versions on the arXiv, we used ChatGPT for proofreading before final submission to a journal. 
\section{Preliminaries}

Wigner's semicircle distribution $\rhosc$ and its Stieltjes transform $\msc (z)$ are defined by
\beq \label{eqn:sc-def}
\rhosc (x) := \1_{ \{ |x| \leq 2 \} } \frac{ \sqrt{4-x^2}}{2 \pi }, \quad \msc (z) := \int \frac{ \rhosc (x)}{x -z } \d x = \frac{ -z  + \sqrt{z^2-4}}{2} ,
\eeq
which describe the asymptotic empirical eigenvalue distribution of the G$\beta$E \cite{erdHos2013local}.
We make note of the particular value $\msc(2) = -1$ which, after re-scaling,  corresponds to the $N^{1/3}$ term  on the RHS of \eqref{eqn:intro-1} and \eqref{eqn:yhat-def}. 
We define the quantiles $\gamma_i$ by the equation
\beq
\int_{\gamma_i}^2 \rhosc (x) \d x = \frac{i}{N} , \qquad i \in \{ 0, 1, \dots, N \}
\eeq
It is easy to see that
\beq \label{eqn:quantiles-approx}
\gamma_i = 2 - \left( \frac{ 3 \pi i }{2 N } \right)^{2/3} \left(1 + \O \left( \frac{ i^{2/3}}{N^{2/3}} \right) \right)
\eeq
for $i \leq c N$ for some $c>0$. In light of the convergence of the largest eigenvalues of the G$\beta$ E to the Airy$_\beta$ point process, we introduce the re-scaled quantity
\beq \label{eqn:rho-def}
\rho (E) := N^{1/3} \rhosc(2- E N^{-2/3} ) = \frac{\sqrt{E}}{\pi} \left(1 + \O ( E N^{-2/3} ) \right) ,
\eeq
the asymptotic holding for $0 \leq E \leq c N^{2/3}$, some $c>0$. We also define the quantities
\beq \label{eqn:hatgam-def}
\tilgam_j :=N^{2/3} (2 - \gamma_j), \quad 0 \leq j \leq N , \qquad \hatgam_i :=\left( \frac{3 \pi i }{2} \right)^{2/3} , \quad i \geq 0.
\eeq
The following two results are proven in Appendices \ref{sec:rigi-proof} and \ref{sec:iso}, respectively. The first estimate is a relatively standard consequence of the optimal local law of Bourgade, Mody and Pain \cite{bourgade2022optimal}. One could attempt to quantify the dependence of the constant on $p$, but we make no attempt to do so. More importantly, the $| \log i|$ dependence is not expected to be optimal; the optimal result should be $| \log i |^{1/2}$. 
\bep \label{prop:rigidity}
Let $x_i$ be the eigenvalues of the re-scaled Gaussian $\beta$-ensemble as in \eqref{eqn:GBE-rescaled}. For any $p \in \nn$ there is a constant $C_p >0$ so that
\beq \label{eqn:rigid-moment}
\left( \ee \left| x_i - \hatgam_i \right|^p \right)^{1/p} \leq C_p \frac{ 1+|\log i|}{ i^{1/3}}
\eeq
uniformly for $1 \leq i \leq N^{1/10}$. For any $\eps >0$ and $ D>0$ we have
\beq \label{eqn:usual-rigid}
\pp \left[ |x_i - \tilgam_i | > \frac{ N^{\eps}}{ \min\{ i^{1/3} , (N+1-i )^{1/3} \} } \right] \leq N^{-D} .
\eeq
\eep

\bel \label{lem:iso}
Suppose that $x_1 < x_2 < \dots < x_N$ are points obeying the estimates \eqref{eqn:usual-rigid} and let $\{z_i\}_i$ be i.i.d. Gamma random variables with distribution as in Section \ref{sec:def-processes}, independent of the $\{x_i \}_{i=1}^N$. Then for all $\eps >0$ and $D>0$ we have
\beq
\pp \left[ \left| \sum_{i=1}^N \frac{z_i}{ E - x_i} - N^{1/3} \msc (2  - EN^{-2/3} ) \right| > \frac{N^{\eps}}{ |E|^{1/4}} \right] \leq N^{-D} 
\eeq
for all $- N^{2/3} \leq E \leq - N^{\eps}$. 
\eel
From \eqref{eqn:rigid-moment}, the convergence of the Gaussian $\beta$-ensemble to the Airy$_\beta$ point process, and Fatou's lemma we immediately obtain the following.
\bec \label{cor:airy-rig}
Let $\mfa_i$ denote the Airy$_\beta$ point process. For all $p \geq 1$ there is a $C_p$ so that
\beq
\left( \ee \left| \mfa_i - \hatgam_i \right|^p\right)^{1/p}  \leq C_p \frac{1+| \log i | }{ i^{1/3}}
\eeq
\eec

\subsection{Notation}

We say that an event $\F$ holds with overwhelming probability if for all $D>0$ we have $\pp[ \F] \geq 1  - N^{-D}$ for $N$ large enough. For two non-negative sequences $a_N$ and $b_N$ we write $a_N \lesssim b_N$ if there is a constant $C>0$ so that $a_N \leq C b_N$. Sometimes the $a_N$ and $b_N$ may depend on some auxiliary parameter. In this case we write $a_N \lesssim b_N$ if the constant $C>0$ is independent of that parameter. We write $a_N \asymp b_N$ if $a_N \lesssim b_N \lesssim a_N$.

\section{Limiting process} \label{sec:limiting-process}

\subsection{Existence} \label{sec:existence}

In this section we show that the limiting processes exist, which essentially amounts to showing that the function $F(E; \mfa , \bz)$ in \eqref{eqn:F-def} is well-defined for the processes we consider. Throughout this section and the next we will denote $F(E ) = F(E ; \mfa , \bz )$ when there is no cause for confusion.
\bed \label{def:airy}
We say that a simple point process $\mfa_1 < \mfa_2 < \dots $ is Airy-like if there is a $\hat{K}\geq 1$ so that for all $p \geq 1$ we have for some constant $C_p >0$
\beq \label{eqn:airy-def-moment}
\ee[ | \mfa_i - \hatgam_i |^p] \leq C_p \left( \frac{1+| \log i | }{ i^{1/3}} \right)^p
\eeq
for all $i \geq \hat{K}$. 
\eed
Note that we are assuming that $\mfa_1 \in \rr$ so there is some inherent tightness to the first particle (i.e., $\mfa_1 \neq - \infty$ a.s.). We first check that $F$ is well defined in our setting.

\bel \label{lem:series-exists}
Let $\mfa = ( \mfa_1, \dots )$ be Airy-like. Let $\bz = ( z_1, \dots )$ be as above. Then there is an event with probability $1$ on which the limits
\beq
\lim_{n \to \infty} \left( \sum_{i=1}^n \frac{z_i}{ E - \mfa_i } + \int_0^{\hatgam_n} \frac{1}{ \pi \sqrt{x} } \d x \right)
\eeq
exist simultaneously 
for all $E \in \rr \backslash \{ \mfa \}$. 
\eel
\proof From the explicit density of the $z_i$, the definition of $\hatgam$ in \eqref{eqn:hatgam-def}, the Airy-like assumption \eqref{eqn:airy-def-moment}, and Borel-Cantelli it is easy to see that there is a $c_1 >0$ so that almost surely there is an $L \geq 1$ so that
\beq \label{eqn:as-event}
|z_i | \leq i^{1/10} , \qquad  | \mfa_i - \mfa_1 | \geq c_1 i^{2/3}
\eeq
for all $i \geq L$. In particular, we used that
\beq
\pp\left[ | \mfa_i - \hatgam_i | > \frac{1}{2} | \hatgam_i| \right] \lesssim \frac{1}{ (i^{2/3} )^p} \ee\left[ | \mfa_i - \hatgam_i |^p \right] \lesssim \frac{ | \log i |^p}{ i^p}
\eeq
for any $p \geq 1$ and $i$ sufficiently large. On the event of \eqref{eqn:as-event} one finds that the series
\beq \label{eqn:aaa-2}
\sum_{i=1}^\infty z_i \left| \frac{1}{ E - \mfa_i } - \frac{1}{\mfa_1 - \mfa_i -1 } \right| \leq \sum_{i=1}^\infty z_i \frac{1+|E|+|\mfa_1|}{ |E-\mfa_i| |\mfa_1 - \mfa_i -1 | }
\eeq
is absolutely convergent simultaneously for all $E \in \rr \backslash \{ \mfa \}$, because the denominator on the RHS grows like $i^{4/3}$ for all sufficiently large $i$. It therefore suffices to prove convergence of
\beq
\lim_{n \to \infty} \left( \sum_{i=1}^n \frac{z_i}{ \mfa_1 - \mfa_i-1 } + \int_0^{\hatgam_n} \frac{1}{ \pi \sqrt{x} } \d x \right)
\eeq
almost surely. With $c_1>0$ as above, we first show convergence (as $n \to \infty$) of,
\begin{align}
& \sum_{i=1}^n \left( \frac{z_i}{ \mfa_1 - \mfa_i -1 }  + z_i \int_{\hatgam_{i-1}}^{\hatgam_i} \frac{1}{ \pi \sqrt{x}} \d x \right) \\
= &  \sum_{i=1}^n \1_{\{  c_1 i^{2/3} > | \mfa_i - \mfa_1 | \} } \left( \frac{z_i}{ \mfa_1 - \mfa_i -1 }  + z_i \int_{\hatgam_{i-1}}^{\hatgam_i} \frac{1}{ \pi \sqrt{x}} \d x  \right) \notag\\
+ & \sum_{i=1}^n \1_{ \{ c_1 i^{ 2/3} \leq  | \mfa_i - \mfa_1 | \} } \left( \frac{z_i}{ \mfa_1 - \mfa_i -1 } + z_i \int_{\hatgam_{i-1}}^{\hatgam_i} \frac{1}{ \pi \sqrt{x}} \d x  \right) 
\end{align} 
The limit of the second line always exists as there are only finitely many non-zero terms in the series, almost surely. We now show that the sum in the last line is absolutely convergent, almost surely. By the definition of the quantiles we have
\begin{align} \label{eqn:aaa-1}
&\sum_{i=1}^n \1_{ \{c_1 i^{ 2/3} \leq  | \mfa_i - \mfa_1 |\} } \left| \frac{z_i}{ \mfa_1 - \mfa_i -1 } + z_i \int_{\hatgam_{i-1}}^{\hatgam_i} \frac{1}{ \pi \sqrt{x}} \d x \right| \notag\\
\lesssim &  \sum_{i=1}^n \1_{ \{c_1 i^{ 2/3} \leq  | \mfa_i - \mfa_1 |\} }  z_i \int_{\hatgam_{i-1}}^{\hatgam_i} \frac{ |x- \mfa_i| + | \mfa_1| +1}{x|\mfa_1 - \mfa_i -1 |} \sqrt{x} \d x 
\lesssim  \sum_{i=1}^\infty {z_i} \frac{1 + |\mfa_1| + | \mfa_i - \hatgam_i|}{i^{4/3}}.
\end{align}
We also used $| \hatgam_i  - \hatgam_{i-1} | \lesssim i^{-1/3} \leq 1$.
By the Airy-like property we have (here $\hat{K}$ is the constant from Definition \ref{def:airy})
\beq \label{eqn:airy-finite-1}
\sum_{i=\hat{K}}^\infty \ee\left[ |z_i| \frac{ 1+ | \mfa_i - \hatgam_i|  }{ i^{4/3}} \right] \lesssim 1 .
\eeq
From this and the fact that $|\mfa_1| < \infty$ a.s. we have that the RHS of \eqref{eqn:aaa-1} is a.s. absolutely convergent.  Finally, we see that
\beq
\lim_{n \to \infty}  \sum_{i=1}^n ({ 1 - z_i}) \int_{\hatgam_{i-1}}^{\hatgam_i} \frac{1}{ \pi \sqrt{x}} \d x 
\eeq
converges by the Kolmogorov Three Series test since $\ee[z_i] =1$, and the coefficient  of $(1-z_i)$ in the series is $ \O ( i^{-2/3} )$ and so is square-summable. \qed

We now show that $F$ is sufficiently regular so that we can construct the process $\mfb$ as discussed in Section \ref{sec:def-processes}. 
\bel \label{lem:limits} Let $\mfa$ and $\bz$ be as above. There is an event with probability $1$ on which the function $E \to F (E ; \mfa , \bz )$ is  continuous and strictly decreasing on the connected components of $\rr \backslash \{ \mfa \}$ and
\beq \label{eqn:limits}
\lim_{E \to \mfa_i^-} F (E ; \mfa , \bz ) = - \infty , \qquad \lim_{E \to \mfa_{i-1}^+} F (E ; \mfa, \bz ) = + \infty .
\eeq
\eel
\proof The continuity and strict monotonicity follow from taking differences $F(E_1) - F(E_2)$ and checking that the series
\beq
\sum_{j=1}^\infty \frac{z_j}{ (E_1 - \mfa_j) (E_2 - \mfa_j) }
\eeq
converges absolutely for $E_1, E_2 \in ( \mfa_{i-1}, \mfa_i )$ on the almost-sure event that \eqref{eqn:as-event} holds for all large $j$. In particular, every term in the series is positive. 

It remains to compute the limits in \eqref{eqn:limits}. First consider $E \in ( \mfa_i , \mfa_{i+1} )$ for some $i \geq 1$. We write,
\beq \label{eqn:limits-1}
F(E) - F ( \mfa_1-1 ) = \sum_{j=1}^{i+1} z_j \left( \frac{1}{ E- \mfa_j } - \frac{1}{ \mfa_1 - \mfa_j -1 } \right) +  \sum_{j=i+2}^\infty  z_j \left( \frac{1}{ E- \mfa_j } - \frac{1}{ \mfa_1 - \mfa_j -1 } \right)
\eeq
Since we know that almost surely there is an $L >0$ so that $|\mfa_j - \mfa_1 | \geq c_1 j^{2/3}$ and $|z_j| \leq j^{1/10}$ for all $j \geq L$, it is easy to check that the second term on the RHS is uniformly bounded for $E \in (\mfa_{i} , \mfa_{i+1} )$, a.s. (the bound is random, but holds uniformly in $E \in (\mfa_i, \mfa_{i+1})$). This argument is similar to our previous estimation of \eqref{eqn:aaa-2} and also uses the fact that $\mfa$ is simple so that $\mfa_{i+2}  > \mfa_{i+1}$.  Therefore we can compute the limits $E \to \mfa_i^+$ and $E \to \mfa_{i+1}^-$ from the first term on the RHS of \eqref{eqn:limits-1}, finding they are $+\infty$ and $- \infty$ respectively (recall $\bz_j >0$ for all $j$). The limit $E \to \mfa_1^-$ is similar. It remains to show that
\beq
\lim_{E \to - \infty} F(E) = \infty ,
\eeq
almost surely. Since $F(E)$ is a monotonic function it suffices to prove that for all $\eps >0$ there is an $E \geq 1$ so that
\beq \label{eqn:F-bd-aa-1}
\pp\left[ F(\mfa_1 - E) > ( \eps )^{-1} \right] \geq 1- \eps .
\eeq
We write for $E \geq 1$,
\begin{align} \label{eqn:F-bd-aa-2}
F( \mfa_1 - E) &= \left( \sum_{i=1}^\infty \frac{z_i -1 }{\mfa_1 - \mfa_i - E } \right) + \left( \sum_{i=1}^\infty \frac{1}{ \mfa_1 - \mfa_i - E} - \int_{\hatgam_{i-1}}^{\hatgam_i} \frac{ \sqrt{x}}{ \pi (-x-E) } \d x \right) \notag\\
&+ \int_0^\infty \left( \frac{1}{x} - \frac{1}{x + E } \right) \frac{ \sqrt{x}}{ \pi } \d x .
\end{align}
From Lemma \ref{lem:calc-1} we see that there is a constant $c>0$ so that the last line of \eqref{eqn:F-bd-aa-2} is larger than $c E^{1/2}$. Call the first two random variables on the RHS of \eqref{eqn:F-bd-aa-2} $X(E)$ and $Y(E)$ respectively. In order to prove \eqref{eqn:F-bd-aa-1} it therefore suffices to prove that for all $\eps >0$ there is a constant $C_\eps$ so that
\beq \label{eqn:XY-bd}
\sup_{E \geq 1 } \pp\left[ |X(E)| + |Y(E)| > C_\eps \right] \leq \eps .
\eeq
We now turn to the proof of this estimate, which will complete the proof of the lemma. We start with $X(E)$. Let $c_1 >0$ be the constant from \eqref{eqn:as-event}. We write
\begin{align}
X(E) = \left( \sum_{i=1}^\infty \1_{ \{ | \mfa_i - \mfa_1 | \geq c_1 i^{2/3} \} } \frac{z_i -1 }{\mfa_1 - \mfa_i - E }  \right) + \left( \sum_{i=1}^\infty \1_{ \{ | \mfa_i - \mfa_1 | < c_1 i^{2/3} \} } \frac{z_i -1 }{\mfa_1 - \mfa_i - E } \right)
\end{align}
Since $E\geq 1$, the second term is bounded above by
\beq
\left|  \sum_{i=1}^\infty \1_{ \{ | \mfa_i - \mfa_1 | < c_1 i^{2/3} \} } \frac{z_i -1 }{\mfa_1 - \mfa_i - E } \right| \leq \sum_{i=1}^\infty \1_{ \{ | \mfa_i - \mfa_1 | < c_1 i^{2/3} \} } ( |z_i| +1 )
\eeq
which is a.s. a finite random variable (as the series has only finitely many terms, a.s.). By the independence of the $z_i$ from the $\mfa_i$ we see for the first term (by taking the expectation over the $z_i$ first) that
\begin{align}
\ee \left[\left( \sum_{i=1}^\infty \1_{ \{ | \mfa_i - \mfa_1 | \geq c_1 i^{2/3} \} } \frac{z_i -1 }{\mfa_1 - \mfa_i - E }  \right)^2  \right] \lesssim \sum_{i=1}^\infty \frac{1}{ i^{4/3}} \lesssim 1
\end{align}
uniformly for $E \geq 1$. This completes the bounding of $X(E)$. We turn to $Y(E)$. We write it as,
\begin{align}
  Y(E)  =&    \sum_{i=1}^\infty \1_{ \{ |\mfa_i - \mfa_1| < c_1 i^{2/3} \} } \left(\frac{1}{ \mfa_1 - \mfa_i - E} - \int_{\hatgam_{i-1}}^{\hatgam_i} \frac{ \sqrt{x}}{ \pi (-x-E) } \d x \right)  \notag\\
 + &   \sum_{i=1}^\infty \1_{ \{ |\mfa_i - \mfa_1| \geq c_1 i^{2/3} \} } \left(\frac{1}{ \mfa_1 - \mfa_i - E} - \int_{\hatgam_{i-1}}^{\hatgam_i} \frac{ \sqrt{x}}{ \pi (-x-E) } \d x \right)
\end{align}
By a similar argument to the $X(E)$ case, the first term on the RHS is bounded above by an $E$-independent random variable that is almost surely finite. By an argument similar to \eqref{eqn:aaa-1}, the second term is bounded above by,
\beq
\sum_{i=1}^\infty \1_{ \{ |\mfa_i - \mfa_1| \geq c_1 i^{2/3} \} } \left|\frac{1}{ \mfa_1 - \mfa_i - E} - \int_{\hatgam_{i-1}}^{\hatgam_i} \frac{ \sqrt{x}}{ \pi (-x-E) } \d x \right| \lesssim \sum_{i=1}^\infty \frac{ | \mfa_1| + | \mfa_i - \hatgam_i | + 1 }{i^{4/3}} ,
\eeq
uniformly in $E \geq 1$. The RHS is a.s. finite by an estimate similar to \eqref{eqn:airy-finite-1}. This completes the proof of \eqref{eqn:XY-bd} and so the proof of the lemma is complete. \qed

From Lemmas \ref{lem:series-exists} and \ref{lem:limits} we see that if $\mfa$ is Airy-like and $\bz$ is as above, then $(\mfa, \bz)$ is admissible as defined in Section \ref{sec:def-processes}. Therefore, the particle process $\mfb = \Phi (\mfa, \bz, \alpha)$ is well-defined. Proposition \ref{prop:intro-well-defined} is then an immediate corollary of the following.
\bel
If $\mfa$ is Airy-like and $\bz$ is as above, then the  particle process $\mfb$ defined by $\mfb = \Phi ( \mfa , \bz , \alpha )$ is well-defined and interlaces $\mfa$. Moreover, it is Airy-like. 
\eel
\proof The estimate \eqref{eqn:airy-def-moment} for the $\mfb_i$ follows from the interlacing property and the same estimate for the $\mfa_i$, after increasing $\hat{K}$ by $1$ and increasing the constants $C_p$ appropriately. Note we also use the fact that $| \hatgam_i - \hatgam_{i-1} | \lesssim (i)^{-1/3}$. \qed

\subsection{Approximation}

Let $\mfa$ be an Airy-like point process and $\bz$ as above. We consider two processes, $\mfb = \Phi ( \mfa, \bz, \alpha )$ as above and for every $m$, the process $\hfb = \{ \hfb_i \}_{i=1}^{m}$ defined as  the solutions of
\beq \label{eqn:airy-finite}
\sum_{i=1}^m \frac{ z_i}{ E - \mfa_i} + \int_0^{ \hatgam_m} \frac{1}{ \pi \sqrt{x}} \d x = \alpha ,
\eeq
that satisfy $\hfb_1 < \mfa_1 < \dots < \hfb_m <\mfa_m $. 
We tacitly assume that $m > 100 | \alpha |^3$ so that $\hfb_1$ is well defined (i.e., this ensures that the deterministic term on the LHS is strictly larger than $\alpha$). Note we suppress the dependence of $\hfb$ on $m$ in the notation. We first prove a certain tightness estimate for the first particle $\hfb_1$. 
\bel \label{lem:approx-tight}
Let $\eps >0$. There exist $M, m_0 >0$ so that for all $m \geq m_0$ we have
\beq
\pp[ | \hfb_1 | > M ] \leq \eps .
\eeq
\eel
\proof It suffices to prove that for all $\eps >0$, there exist $m_0$ and $M \geq 1$ so that for all $m \geq m_0$
\beq \label{eqn:finite-tight-a2}
\pp\left[ \sum_{i=1}^m \frac{z_i}{\mfa_1 - \mfa_i -M} + \int_0^{\hatgam_m} \frac{1}{ \pi \sqrt{x}} \d x > \alpha \right] \geq 1 - \eps ,
\eeq
because on the event inside the probability on the LHS we have $\mfa_1 - M < \hfb_1 < \mfa_1$. 
We first write
\begin{align} \label{eqn:finite-tight-a1}
\sum_{i=1}^m \frac{z_i}{\mfa_1 - \mfa_i -M}  & = \sum_{i=1}^m \left( \frac{z_i}{\mfa_1 - \mfa_i -M}  - z_i\int_{\hatgam_{i-1}}^{\hatgam_i} \frac{\sqrt{x}}{ \pi (-x-M) } \d x \right) \notag\\
&+ \left(\sum_{i=1}^m (z_i-1)\int_{\hatgam_{i-1}}^{\hatgam_i} \frac{\sqrt{x}}{ \pi (-x-M) } \d x \right) + \int_0^{\hatgam_m} \frac{ \sqrt{x}}{ \pi(-x- M)} \d x .
\end{align}
We now bound the first two terms on the RHS. Let $\eps >0$. As seen earlier in the paper there is a $c_1 >0$ so that we may choose $L \geq \hat{K}$ so that with probability at least $1 - \eps$ we have for all $i \geq L$ that $| \mfa_i - \mfa_1 | \geq c_1 i^{2/3}$. Choose $C_1 \geq 1$ so that with probability at least $1- \eps$ we have $|\mfa_1| \leq C_1$ and $\sum_{i=1}^L |z_i| \leq C_1$. Call the event that these estimates hold $\A$. On the event $\A$ we have, using that $M \geq 1$ and by possibly increasing $C_1 >0$ that,
\begin{align}
\left| \sum_{i=1}^m \left( \frac{z_i}{\mfa_1 - \mfa_i -M}  - z_i\int_{\hatgam_{i-1}}^{\hatgam_i} \frac{\sqrt{x}}{ \pi (-x-M) } \d x \right) \right| \leq C_1 \left( 1 +  \sum_{i=\hat{K}}^m z_i\frac{ 1 + | \mfa_i - \hatgam_i |}{  i^{4/3} } \right).
\end{align}
Note that the RHS does not depend on $M$. By an argument similar to \eqref{eqn:airy-finite-1}, the expectation of the second term is bounded above independently of $m$. We conclude that for all $\eps >0$, there is a constant $C_\eps >0$ so that for all $m \geq 1$,
\beq \label{eqn:finite-a3}
\pp \left[ \sup_{ M \geq 1} \left| \sum_{i=1}^m \left( \frac{z_i}{\mfa_1 - \mfa_i -M}  - z_i\int_{\hatgam_{i-1}}^{\hatgam_i} \frac{\sqrt{x}}{ \pi (-x-M) } \d x \right) \right| > C_\eps \right] \leq \eps .
\eeq
It is easy to see that the first term on the second line of \eqref{eqn:finite-tight-a1} is centered and its variance is bounded by a constant independent of $M$. Therefore we see that for all $\eps >0$, there is a constant $C_\eps$ so that for all $M \geq 1$ and $m \geq 1$,
\beq \label{eqn:finite-a4}
\pp\left[ \left| \sum_{i=1}^m (z_i-1)\int_{\hatgam_{i-1}}^{\hatgam_i} \frac{\sqrt{x}}{ \pi (-x-M) } \d x \right| > C_\eps \right] \leq \eps .
\eeq
From the decomposition \eqref{eqn:finite-tight-a1} and the estimates \eqref{eqn:finite-a3} and \eqref{eqn:finite-a4} we conclude that the following statement holds: for all $\eps >0$, there exists a constant $C_\eps >0$ so that for every $M \geq 1$ and $m \geq 1$, there is an event $\A (M, m)$ satisfying $\pp[ \A(M, m)  ] \geq 1 - \eps$ on which the random variable in the probability in the LHS of \eqref{eqn:finite-tight-a2} is bounded below by
\beq
\frac{1}{\pi}\int_0^{\hatgam_m} \left( \frac{1}{x} - \frac{1}{ x + M} \right) \sqrt{x} \d x - C_\eps .
\eeq
If we take now $M \geq C (1+ (C_\eps)^2 + \alpha^2)$ for some large $C>0$, by Lemma \ref{lem:calc-1} the first term above is larger than $C_\eps + |\alpha| +1$ for $m \geq CM^{3/2}$. This establishes the estimate \eqref{eqn:finite-tight-a2} and so completes the proof of the lemma. \qed

The following is a preliminary lemma towards showing that the $\hfb_i$ approximate the $\mfb_i$. 
\bel \label{lem:approx-1}
Let $J \geq 1$ and $\eps >0$. There exists $m_0 \geq 1$ so that for all  $m \geq m_0$ we have
\beq
\pp\left[ \sup_{1 \leq j \leq J} \left| \lim_{n \to \infty} \sum_{i=m}^n \frac{ z_i}{ \mfb_j - \mfa_i} + \int_{\hatgam_{m-1}}^{ \hatgam_n} \frac{1}{ \pi \sqrt{x}} \d x \right| > \eps \right] \leq \eps .
\eeq
\eel
\proof By adjusting the value of $\eps$ it suffices to prove the inequality for a fixed $\mfb_j$. Let us write the quantity in question as
\begin{align} \label{eqn:approx-a1}
 & \sum_{i=m}^n \frac{ z_i}{ \mfb_j - \mfa_i} + \int_{\hatgam_{m-1}}^{ \hatgam_n} \frac{1}{ \pi \sqrt{x}} \d x = \left( \sum_{i=m}^n \left( \frac{ z_i}{ \mfb_j - \mfa_i} - z_i \int_{\hatgam_{i-1}}^{\hatgam_i} \frac{\sqrt{x}}{\pi (-x-1) } \d x \right) \right) \notag\\
+ & \left( \sum_{i=m}^n (z_i -1 )\int_{\hatgam_{i-1}}^{\hatgam_i} \frac{\sqrt{x}}{\pi(-x-1) } \d x  \right) + \left( \int_{\hatgam_{m-1}}^{\hatgam_n} \frac{\sqrt{x}}{\pi} \left( \frac{1}{x} - \frac{1}{x+1} \right ) \d x\right)
\end{align}
It is clear that the last  term on the RHS (which is deterministic) is $\O ( m^{-1/3} )$ uniformly in $n$ so we can ignore it. By the Kolmogorov Three Series test the first term on the second line converges as $n \to \infty$ and by Fatou's lemma,
\begin{align}
\ee \left[ \left( \sum_{i=m}^\infty (z_i -1 )\int_{\hatgam_{i-1}}^{\hatgam_i} \frac{\sqrt{x}}{(-x-1) } \d x  \right)^2 \right] & \leq \liminf_{n \to \infty} \ee \left[ \left( \sum_{i=m}^n (z_i -1 )\int_{\hatgam_{i-1}}^{\hatgam_i} \frac{\sqrt{x}}{(-x-1) } \d x \right)^2 \right] \notag\\
 & \lesssim m^{-1/3}
\end{align}
by direct calculation. Therefore, for all $\eps >0$, there exists an $m_0$ so that for all $m \geq m_0$,
\beq
\pp\left[ \left| \sum_{i=m}^\infty (z_i -1 )\int_{\hatgam_{i-1}}^{\hatgam_i} \frac{\sqrt{x}}{\pi (-x-1) } \d x  \right| > \eps \right] \leq \eps .
\eeq
We now turn to the first term on the RHS of \eqref{eqn:approx-a1}.  From the Airy-like property (recall Definition \ref{def:airy}) and the fact that $\mfb_j \in \rr$ a.s., we see that there exists a $c_1>0$ so that for all $\eps >0$, there is an  $L \geq \hat{K}$ so that for all $i \geq L$ we have $\mfa_i \geq c_1 i^{2/3}$ and $|\mfb_j| \leq L$ with probability at least $1 - \eps$. On the event this holds we have for all $m \geq m_0  =m_0 (L)$
\begin{align}
 \sum_{i=m}^\infty \left| \frac{ z_i}{ \mfb_j - \mfa_i} - \frac{z_i}{ \pi} \int_{\hatgam_{i-1}}^{\hatgam_i} \frac{\sqrt{x}}{(-x-1) } \d x \right| \leq C_\eps \sum_{i=m_0}^\infty z_i \frac{ L +1 + | \mfa_i - \hatgam_i | }{i^{4/3}} ,
\end{align}
for a constant $C_\eps$ depending on $\eps >0$ through our choice of $L$. 
The expectation of the series on the RHS is $\O ( (m_0)^{-1/3} )$. We conclude the desired inequality by Markov's inequality, after taking $m_0$ possibly larger. \qed

We now can show that the $\hfb_j$ are  a good approximation of the $\mfb_j$ with high probability. 
\bep \label{prop:approx-1}
Let $\eps >0$ and $J \geq 1$. There exists an $m_0 \geq 0$ so that for all $m \geq m_0$ we have
\beq
\pp\left[ \sup_{1 \leq j \leq J} | \mfb_j  - \hfb_j | > \eps \right] \leq \eps .
\eeq
\eep
\proof By adjusting the value of $\eps$ it suffices to prove the inequality for a fixed $j$. Let $\eps >0$ be fixed throughout the proof. If $j \geq 2$, we see that there is a $\delta >0$ so that
\beq \label{eqn:approx-a2}
\frac{ z_j}{ (\mfa_j - \mfa_{j-1} )^2 } \geq \delta
\eeq
with probability at least $1- \eps$. If $j=1$ then by the fact that $|\mfb_1| < \infty$ a.s. and Lemma \ref{lem:approx-tight} we have that there exist  $m_0$ and  $\delta >0$ so that
\beq \label{eqn:approx-a3}
\frac{z_1}{ ( \mfa_1 - \mfb_1 ) ( \mfa_1 - \hfb_1 ) } \geq \delta
\eeq
with probability at least $1- \eps$ for all $m \geq m_0$. Define $\A (m)$ to be the event that  \eqref{eqn:approx-a2} or \eqref{eqn:approx-a3} holds in the case $j\geq 2$ or $j=1$, respectively.  By Lemma \ref{lem:approx-1} we see that there is an $m_1 \geq 1$ so that for all $ m \geq m_1$ there is an event $\B (m)$ on which
\beq
\lim_{n \to \infty} \left|  \sum_{i=m+1}^n \frac{ z_i}{ \mfb_j - \mfa_i} + \int_{\hatgam_{m}}^{ \hatgam_n} \frac{1}{ \pi \sqrt{x}} \d x \right| < \eps \delta 
\eeq
and $\pp[ \B (m)] \geq 1- \eps$. By taking the difference between the equations satisfied by $\mfb_j$ and $\hfb_j$ we find
\beq
| \mfb_j - \hfb_j | \sum_{i=1}^m \frac{z_i}{( \mfa_i - \mfb_j ) ( \mfa_i - \hfb_j ) } \leq \eps \delta
\eeq
on the event $\B (m)$.  Here we used that every term in the series on the LHS is positive. Using this positivity we see that on the event $\A(m)$ the series on the LHS is bounded below by $\delta$. This yields the claim. \qed

\section{Pre-limiting process}

In this section we consider various properties of the sequence of solutions to \eqref{eqn:yhat-def}. The following encapsulates the useful properties that the pre-limiting processes we consider have.

\bed Let $\ell \geq 0$ be fixed. For every $N \geq 1$ we consider a  point  process of $N +\ell$ particles that we denote by $x_1 < x_2 < \dots < x_{N + \ell}$.\footnote{Technically, we should denote the particles by $x_i^{(N)}$. That is, we should include a superscript to indicate the dependence on the index $N$. We suppress this dependence for notational simplicity.} We say they are rigid if the following holds. For all $\eps >0$ and $D>0$ we have
\beq \label{eqn:rig-def-1}
\pp\left[ | x_i - \tilgam_i | > \frac{ N^{\eps}}{ \min \{ i^{1/3} , (N+1 -i )^{1/3} \} } \right] \leq N^{-D}
\eeq
for all $1 \leq i \leq N$ and $N$ large enough. Furthermore, there is a $\hat{K} \geq 1$ so that for all $p \geq 1$ there is a $C_p$ so that for all $\hat{K} \leq i \leq N^{1/10}$
\beq \label{eqn:rig-def-2}
\ee \left[ | x_i - \hatgam_i |^p \right] \leq C_p \left( \frac{1+ |\log i| }{i^{1/3}} \right)^p ,
\eeq
and $x_1$ is a tight random variable.
\eed

Fix now a rigid process $x_i$, a vector $\bz$ as earlier, a centered Gaussian random variable $g$, and $\alpha \in \rr$. We will consider the following two point processes, $\{ y_i \}_{i=1}^{N+\ell+1}$ and $\{ \tily_i \}_{i=1}^{N+\ell}$. First the $y_i$ are the $N+\ell+1$ unique solutions to
\beq \label{eqn:y-def-aa}
\alpha = \sum_{i=1}^{N+\ell} \frac{z_i}{ y -x_i} + N^{1/3} - (y N^{-1/3} + g N^{-1/6} ),
\eeq
in increasing order s.t. $y_1 < x_1 < y_2 < x_2 < \dots < x_{N+\ell} < y_{N+\ell+1}$. 
Secondly, $\tily_i$ are the $N+\ell$ unique solutions to
\beq \label{eqn:y-def}
\alpha = \sum_{i=1}^{N+\ell} \frac{z_i}{ \tily -x_i} + \int_{0}^{\hatgam_{N^{1/10}}} \frac{1}{ \pi x^{1/2}} \d x + \int_{\tilgam_{N^{1/10}}}^\infty \frac{\rho (x)}{x} \d x .
\eeq
We remark that the $\tily_i$  are only defined for a deterministic $N$ large enough depending only on $\alpha$ (or else there may not be a solution $\tily_1 < x_1$ of the above equation). The $\tily_i$ are an intermediate process that removes the effect of the $y$ and $g$ on the RHS of \eqref{eqn:y-def-aa} and approximates $\rhosc$ by the square root function close to the spectral edge.     Furthermore by \eqref{eqn:quantiles-approx} and \eqref{eqn:rho-def} we have,
\beq \label{eqn:determ-edge}
\left| \int_{0}^{\hatgam_{N^{1/10}}} \frac{1}{ \pi x^{1/2}} \d x - \int_0^{\tilgam_{N^{1/10}}} \frac{ \rho (x)} {x} \d x \right| \lesssim N^{-1/10}
\eeq
which also motivates the RHS of \eqref{eqn:y-def}. 
\bel
If the $x_i$ are rigid, both the $y_i$ and $\tily_i$ are also rigid.
\eel
\proof We only do the proof in the case of $y_i$, the $\tily_i$ being similar. Let $\delta \in (0, \frac{1}{10} )$. We first prove with overwhelming probability $y_1 >  - N^{\delta}$. It suffices to show that with overwhelming probability,
\beq
\sum_{i=1}^{N+\ell} \frac{z_i}{ y_* -x_i} + N^{1/3} - (y_* N^{-1/3} + g N^{-1/6} ) > \alpha 
\eeq
for $y_* = - N^{\delta}$. From the formula in \eqref{eqn:sc-def} one finds that
\beq
\msc (2 + E) \geq -1 + c E^{1/2}
\eeq
for some $c>0$ and all $0 < E < 1$. Therefore, from Lemma \ref{lem:iso} we have with overwhelming probability,
\beq
\sum_{i=1}^{N+\ell} \frac{z_i}{ y_* -x_i} = N^{1/3} \msc (2 - N^{-2/3} y_* ) + \O ( N^{-\delta/10} ) \geq - N^{1/3} + c |y_*|^{1/2} - N^{-\delta/11} \geq - N^{1/3} + N^{\delta/3} .
\eeq
We also used that $x_i \gtrsim N^{2/3}$ with overwhelming probability if $i \geq N$ to disregard the particles $\{ x_i \}_{i=N+1}^{N+\ell}$ for which we do not have the estimates \eqref{eqn:rig-def-1}. Therefore, since $-N^{\delta}  <y_1 < x_1$ with overwhelming probability, we deduce that \eqref{eqn:rig-def-1} holds for the case $i=1$. It automatically holds for $2 \leq i \leq N$ by interlacing.

 Similarly, by interlacing we see that since \eqref{eqn:rig-def-2} holds for $x_i$ with $i \geq \hat{K}$, it holds for $y_i$ with $i \geq \hat{K}+1$. It remains to check that $y_1$ is tight. Since $x_1$ is tight it suffices to show that for all $\eps >0$, there is an $M \geq 1$ so that
\beq \label{eqn:yhat-tight}
\sum_{i=1}^{N+\ell} \frac{z_i}{ x_1 -x_i - M} + N^{1/3} - ((x_1 -M) N^{-1/3} + g N^{-1/6} ) > \alpha 
\eeq
with probability at least $1- \eps $, for $N$ sufficiently large depending on $\eps>0$. 
By Lemma \ref{lem:aux-rigidity} and  \eqref{eqn:determ-edge} we have
\begin{align}
\sum_{i=1}^{N+\ell} \frac{ z_i}{ y -x_i } + N^{1/3} =\left(  \sum_{i=1}^{N^{1/10}} \frac{z_i}{ y - x_i} + \int_{0}^{\hatgam_{N^{1/10}}} \frac{1}{ \pi x^{1/2}} \d x \right) + \O ( N^{-1/100} )
\end{align}
uniformly for $|y| \leq N^{1/100}$ with overwhelming probability. Similarly to \eqref{eqn:finite-tight-a1} we decompose
\begin{align}
 &\sum_{i=1}^{N^{1/10}} \frac{z_i}{ x_1 - x_i - M} + \int_{0}^{\hatgam_{N^{1/10}}} \frac{1}{ \pi x^{1/2}} \d x =\left( \sum_{i=1}^{N^{1/10}} \left( \frac{z_i}{x_1 - x_i -M}  - z_i\int_{\hatgam_{i-1}}^{\hatgam_i} \frac{\sqrt{x}}{ \pi (-x-M) } \d x \right) \right) \notag\\
+& \left(\sum_{i=1}^{N^{1/10}} (z_i-1)\int_{\hatgam_{i-1}}^{\hatgam_i} \frac{\sqrt{x}}{ \pi (-x-M) } \d x \right) + \int_0^{\hatgam_{N^{1/10}}} \frac{\sqrt{x}}{\pi} \left( \frac{1}{x} - \frac{1}{ x + M } \right) \d x .
\end{align}
The variance of the first term on the second line is bounded uniformly in $N$, independently of $M \geq 1$. Therefore for all $\eps >0$ there is a $C_\eps >0$ s.t.
\beq
\pp\left[ \left|\sum_{i=1}^{N^{1/10}} (z_i-1)\int_{\hatgam_{i-1}}^{\hatgam_i} \frac{\sqrt{x}}{ \pi (-x-M) } \d x \right| > C_\eps \right] \leq \eps ,
\eeq
uniformly in $N$ and $M \geq 1$. Let $\eps >0$. Using the estimates \eqref{eqn:rig-def-2} and tightness of $x_1$, one sees that there is a $c_1 >0$ and an $L = L_\eps \geq \hat{K}$ so that
\beq
|x_1| \leq L, \qquad |x_i - x_1 | \geq c_1 i^{2/3}
\eeq
with probability at least $1-\eps$, for all $L \leq i \leq N^{1/10}$. Furthermore, there is a $C_1 >0$ so that $\sum_{i=1}^L |z_i| \leq C_1$ with probability at least $1- \eps$.  On this event, we may bound 
\begin{align}
\left| \sum_{i=1}^{N^{1/10}} \frac{z_i}{x_1 - x_i -M}  - z_i\int_{\hatgam_{i-1}}^{\hatgam_i} \frac{\sqrt{x}}{ \pi (-x-M) } \d x \right| \leq C_1 \left( 1 + \sum_{i=\hat{K}}^{N^{1/10}} \frac{ z_i (L + | x_i - \hatgam_i | )}{ i^{4/3}} \right ) 
\end{align}
after possibly increasing $C_1 >0$. Note that $C_1$ does not depend on $M$. The expectation of the RHS is bounded above independently of $N$. Therefore, collecting the above estimates we have shown that for all $\eps >0$, there is a $C_\eps >0$, so that for every $N^{1/100} \geq M \geq 1$ and $N \geq C_\eps$, there is an event $\A (M, N)$ on which
\begin{align}
& \sum_{i=1}^{N+\ell} \frac{z_i}{ x_1 -x_i - M} + N^{1/3} - ((x_1 -M) N^{-1/3} + g N^{-1/6} ) \notag\\
\geq &   - C_\eps + \int_0^{\hatgam_{N^{1/10}}} \frac{\sqrt{x}}{\pi} \left( \frac{1}{x} - \frac{1}{ x + M } \right) \d x.
\end{align}
and $\pp[ \A (M, N)] \geq 1- \eps$. The integral on the RHS is larger than $c M^{1/2}$ for all $1 \leq M \leq N^{1/100}$ by Lemma \ref{lem:calc-1}. This yields the desired estimate \eqref{eqn:yhat-tight} and completes the proof. \qed

The following lemma allows us to pass to the intermediate $\tily_j$.

\bel \label{lem:approx-2}
There is a small $c>0$ so that for any $J$ we have
\beq
|y_j - \tily_j | \leq N^{-c}
\eeq
with probability at least $1 - N^{-c}$ for all $1 \leq j \leq J$. 
\eel
\proof By using rigidity of $\tily_i$, $y_i$, $x_i$, Lemma \ref{lem:aux-rigidity} and taking the difference of the two equations defining the solutions $\tily_j$ and $y_j$ we get,
\beq
\sum_{i=1}^{{ N+\ell}} \frac{ z_i}{ \tily_j - x_i} - \frac{z_i}{ y_j -x_i} = \O ( N^{-1/100} )
\eeq
with overwhelming probability. Therefore
\beq
C N^{-1/100} \geq | \tily_j - y_j| \sum_{i=1}^{{ N+\ell}} \frac{z_i}{ ( \tily_j -x_i) ( y_j - x_i)}\geq \frac{| \tily_j - y_j| z_1}{N^{\eps}}
\eeq
for any $\eps >0$ with overwhelming probability, using the rigidity of the $x_i$ and $y_1, \tily_1$. Since $\pp[ |z_1| \leq N^{- \eps } ] \lesssim N^{-\beta\eps/2}$ the claim follows after choosing $\eps >0$ sufficiently small. \qed 

\subsection{Approximation}

We recall that $x_1 < \dots < x_{N+\ell}$ is a rigid process and the $\tily_i$ are the solutions to \eqref{eqn:y-def}. For any fixed $m$ we define (for $N^{1/10} \geq m)$ the solutions $\{ \haty_i\}_{i=1}^{m}$ to the equation
\beq \label{eqn:tily-def}
\sum_{i=1}^m \frac{z_i}{ \haty_j - x_i } + \int_0^{\hatgam_m} \frac{ 1}{ \pi \sqrt{x} } \d x = \alpha .
\eeq
We tacitly assume that $m$ is sufficiently large depending on $\alpha$ so that $\haty_1$ is well-defined. 

The following is proven in an almost identical manner to Lemma \ref{lem:approx-tight}; we omit the proof.
\bel
Let $\eps >0$. There exist $M>0$, $m_0$ and $N_0$ so that for all $N \geq N_0$ and $m_0 \leq m \leq N^{1/10}$ we have
\beq
\pp\left[ | \haty_1| > M \right] \leq \eps .
\eeq 
\eel

The analog of Lemma \ref{lem:approx-1} is the following.
\bel
Let $J \geq 1$ and $\eps >0$. There exist $N_0 \geq 1$ and $m_0 \geq 1$ so that for all $m_0 \leq m \leq N^{1/10}$  and $N \geq N_0$ we have
\beq
\pp\left[ \sup_{1 \leq j \leq J} \left| \sum_{i=m}^{N+\ell} \frac{z_i}{ \tily_j - x_i } + \int_{\hatgam_{m-1}}^{\hatgam_{N^{1/10}}} \frac{1}{ \pi x^{1/2} } \d x + \int_{\tilgam_{N^{1/10}}}^\infty \frac{ \rho (x) }{x} \d x   \right| > \eps \right] \leq \eps .
\eeq
\eel
\proof By Lemma \ref{lem:aux-rigidity} it suffices to prove the estimate for
\begin{align}
 &\sum_{i=m}^{N^{1/10}} \left( \frac{ z_i}{ \tily_j - x_i} + \int_{\hatgam_{i-1}}^{\hatgam_i} \frac{1}{ \pi x^{1/2}} \d x \right) \notag\\
=& \left( \sum_{i=m}^{N^{1/10}} \left( \frac{ z_i}{ \tily_j - x_i} +z_i \int_{\hatgam_{i-1}}^{\hatgam_i} \frac{1}{ \pi x^{1/2}} \d x \right) \right) + \left( \sum_{i=m}^{N^{1/10}} (1-z_i  ) \int_{\hatgam_{i-1}}^{\hatgam_i} \frac{1}{ \pi x^{1/2}} \d x  \right)
\end{align}
This is the same as Lemma \ref{lem:approx-1} but we provide the details for the reader's convenience. The second term is centered and its variance is bounded by $\O ( m^{-1/3} )$ and so may be dispensed with. Let $\eps >0$. Then there is an $L \geq \hat{K}$ so that $x_i \geq c i^{2/3}$ for all $N^{1/10} \geq i \geq L$ and $|\tily_j| \leq L$, for all $N$ sufficiently large with probability at least $1- \eps$. On this event we have for all $m \geq C (L+1)^2$ that the first term is bounded by,
\beq \label{eqn:ytil-aa1}
\left| \sum_{i=m}^{N^{1/10}} \frac{ z_i}{ \tily_j - x_i} +z_i \int_{\hatgam_{i-1}}^{\hatgam_i} \frac{1}{ \pi x^{1/2}} \d x \right| \leq \sum_{i=m}^{N^{1/10}}  C \frac{z_i (L + |x_i - \hatgam_i | )}{i^{4/3}} .
\eeq 
The expectation of the term on the RHS is bounded by $C_\eps / m^{1/3}$. So by Markov's inequality there is an $m_0 \geq 1$ for which for all $m \geq m_0$ we have that the RHS of \eqref{eqn:ytil-aa1} is bounded above by $\eps$ with probability at least $1- \eps$. This completes the proof. \qed

By a proof similar to that of Proposition \ref{prop:approx-1} we find the following. We omit the details.
\bep \label{prop:approx-3}
Let $\eps >0$ and $J \geq 1$. Then there exist an $m_0 \geq 1$ and $N_0 \geq 1$ so that for all $N \geq N_0$ and $m_0 \leq m \leq N^{1/10}$ we have that
\beq
\pp\left[ \sup_{1 \leq j \leq J } | \haty_j- \tily_j | > \eps \right] \leq \eps .
\eeq
\eep

\section{Proof of Theorem \ref{thm:main}} \label{sec:convergence}

For notational simplicity we just do the case of two levels, the general case following from an induction argument and the Markovian structure of the system. Let $x_1 < x_2 < x_3 < \dots $ be a rigid sequence that converges to some Airy-like point process $\mfa$. Let $\bz$ be an independent vector of Gamma distributed random variables as above. By the independence of $\bz$ from $\bx$ and $\mfa$ we may use the same vector $\bz$ in the equations involving $\bx$ and  $\mfa$. That is, we let $\by$ be the solution associated to $\bx$ and $\bz$ through \eqref{eqn:yhat-def} and let $\mfb$ be the solution associated to $\mfa$ and $\bz$ through the map $\Phi$. It suffices to show for bounded Lipschitz $F$ that
\beq
\lim_{N \to \infty} \ee[ F( x_1, \dots x_n , y_1, \dots, y_n ) ] = \ee[ F (\mfa_1, \dots, \mfa_n , \mfb_1, \dots , \mfb_n  ) ]
\eeq
for $n \in \nn$. Since $F$ is Lipschitz we have by Proposition \ref{prop:approx-1}, Lemma \ref{lem:approx-2} and Proposition \ref{prop:approx-3} that it suffices to show that for all sufficiently large $m \geq n$ that 
\beq
(x_1, \dots, x_m , \haty_1, \dots, \haty_{m}, z_1, \dots, z_m )
\eeq
converges in distribution as $N \to \infty$ to
\beq
( \mfa_1, \dots, \mfa_m, \hfb_1, \dots, \hfb_{m}, z_1, \dots, z_m )
\eeq
where $\haty_i$ and $\hfb_i$ are defined by \eqref{eqn:tily-def} and \eqref{eqn:airy-finite}, respectively, with the $m$ there equal to the $m$ here. However, because the defining equations \eqref{eqn:tily-def} and \eqref{eqn:airy-finite} are identical, we see that there is a function $\Psi : \rr^{2m} \to \rr^{m}$ s.t.
\beq
(\hfb_1, \dots , \hfb_{m} )  = \Psi ( \mfa_1, \dots, \mfa_m, z_1, \dots, z_m ) , \quad (\haty_1, \dots, \haty_{m} ) = \Psi ( x_1, \dots, x_m, z_1, \dots, z_m)
\eeq
By the implicit function theorem, the $j$th component of $\Psi$ is a continuous function on the set
\beq
\{ x_1 < \dots < x_m \} \times (0, \infty)^m .
\eeq
Since $\mfa$ is almost surely simple and the $z_i$ are almost surely positive, we get the desired convergence from the continuous mapping theorem. \qed

\appendix

\section{Estimates on eigenvalue locations}

\subsection{Proof of Proposition \ref{prop:rigidity}} \label{sec:rigi-proof}

 Let us introduce some quantities and prove some preliminary results. Throughout we denote the $L^p$ norm of a random variable by $\| X\|_p^p := \ee[ |X|^p]$. We will denote the $L^p ( \rr , \d x)$ norm of a function by $\| f \|_{L^p}$. Throughout we denote by $\{ x_i \}_{i=1}^N$ the particles of the re-scaled G$\beta$E as in the statement of Proposition \ref{prop:rigidity}.

 Define,
\beq
s_N (z) := \sum_{i=1}^N \frac{1}{x_i -z } , \qquad s(z) = - N^{1/3} \msc (2 - N^{-2/3} z ).
\eeq
By an elementary rescaling of \cite[Proposition 3.5]{bourgade2022optimal} there is a $c_* >0$ so that for all $q >0$ there is a constant $C_q > 0$ so that
\beq \label{eqn:BMP-s}
\| s_N (z) - s(z) \|_q \leq \frac{C_q}{\eta}
\eeq
for $z = E + \i \eta$ for any $-100c_* N^{2/3} \leq E \leq (4+100c_*) N^{2/3}$ and $0 < \eta \leq 100c_* N^{2/3}$.  In terms of $\rho (x)$ defined in \eqref{eqn:rho-def} we have that $s(z) = \int \frac{ \rho (E)}{E-z} \d E$.

\bel \label{lem:rig-proof-1}
Let $10 \leq E_0 \leq c_* N^{2/3}$. Let $f$ be a smooth function s.t. $f = 1 $ for $ x\in [-100, E_0]$, and $f = 0$ for $x < - 100 - (E_0)^{-1/2} $ or $x > E_0 + (E_0)^{-1/2}$. Assume that $\| f' \|_1 \leq C$ and $\| f'' \|_1 \leq C (E_0)^{1/2}$ for some $C>0$. For all $p$ there is a $C_p >0$ so that
\beq
\left\| \sum_i f ( x_i ) - \int f(x) \rho (x) \d x \right\|_p \leq C_p \log (E_0) .
\eeq
\eel
\proof Let $A >0$ and let $\chi : \rr \to \rr$ be a symmetric function s.t. $\chi (y) = 1$ for $|y| \leq A$ and $\chi (y) = 0$ for $|y| >2 A$. We may assume that $|\chi' (y) | \lesssim A^{-1}$. By the Helffer-Sj{\"o}strand formula \cite[(11.15)]{erdHos2017dynamical} we have
\begin{align} \label{eqn:rig-proof-a1}
 & \left| \sum_{i} f ( x_i) - \int f(x) \rho (x) \d x \right| \lesssim \int_{\rr^2} | \chi' (y) (f(x) + \i y f'(x) ) | |s_N (x + \i y ) - s(x + \i y ) | \d x \d y \notag\\
+ & \left| \int_{\rr^2} y f''(x) \chi (y) \Im[ s_N (x+ \i y ) - s(x+ \i y ) ] \d x \d y \right| .
\end{align}
Choose $A = E_0$. 
We start with the first term on the RHS. By Minkowski's integral inequality  and \eqref{eqn:BMP-s} we have
\begin{align}
& \left\| \int_{\rr^2} | \chi' (y) (f(x) + \i y f'(x) ) | |s_N (x + \i y ) - s(x + \i y ) | \d x \d y \right\|_p \notag\\
\leq & \int_{\rr^2} | \chi' (y)| ( |f(x) | + |y f'(x) |) \| s_N (x+ \i y ) - s(x + \i y ) \|_p \d x \d y \notag\\
\lesssim &  \frac{1}{A} \int_{x \in \rr , A < |y| < 2 A } A^{-1} ( |f(x)| + A | f '(x) | ) \d x \d y \lesssim \| f\|_{L^1} A^{-1} + \| f' \|_{L^1} \lesssim 1.
\end{align}
We now turn to the term on the second line of \eqref{eqn:rig-proof-a1}. By symmetry  of $\chi$ we can consider the integral over $y >0$.  We will split the $y$ integration into $y < Y_0$ and $y> Y_0$, for some $Y_0 >0$. We choose $Y_0 = (E_0)^{-1/2} \leq A$. Then, since for any Stieltjes transform the function $y \to y \Im[s_N (x + \i y)]$ is increasing we have,
\begin{align} \label{eqn:Ims-a1}
& \left| \int_{0 < y < Y_0} y f''(x) \chi (y) \Im[ s_N (x+ \i y ) - s(x+ \i y ) ] \d x \d y \right| \notag\\
\leq & \int_{ 0 < y < Y_0} |f''(x) | Y_0 ( \Im[ s_N (x + \i Y_0) ] + \Im[ s(x + \i Y_0) ] ) \d x \d y \notag\\
\lesssim & Y_0^2 \int_\rr |f''(x)| (\Im[ s(x+ \i Y_0) ] + | s_N (x+ \i Y_0) - s(x+ \i Y_0) | ) \d x .
\end{align}
For $|x+\i y| \leq 2 c_* N^{2/3}$ we have \cite[Lemma 4.3]{erdHos2013local}
\beq \label{eqn:Ims}
\Im [s (x + \i y ) ] \lesssim \begin{cases} \sqrt{|x+ \i y |} , & x \geq 0 \\ \frac{ |y|}{ |x + \i y |^{1/2}} & x \leq 0 \end{cases}
\eeq
Therefore for $x$ s.t. $f''(x) \neq 0$ we have $\Im[ s(x + \i Y_0) ] \lesssim (E_0)^{1/2}$. Therefore,
\beq
 Y_0^2 \int_\rr |f''(x)| \Im[ s(x+ \i Y_0) ]  \d x \lesssim (Y_0)^2 (E_0)^{1/2} \| f''\|_{L^1} \lesssim 1. 
\eeq
On the other hand, again by Minkowski and \eqref{eqn:BMP-s} we have
\begin{align}
& \left\| Y_0^2 \int_{\rr} |f''(x) | |s_N(x + \i Y_0) - s(x + \i Y_0) | \d x \right\|_p \leq Y_0^2 \int_{\rr} |f''(x)| \| s_N (x+ \i Y_0) - s(x+ \i Y_0) \|_p \d x  \notag\\
\lesssim & Y_0 \| f''\|_{L^1} \lesssim 1.
\end{align}
The last two inequalities show that the last line of \eqref{eqn:Ims-a1} is $\O (1)$. Therefore, we only have left the part of the integration in the second line of \eqref{eqn:rig-proof-a1} where $y > Y_0$. By integration by parts in $x$ and the Cauchy-Riemann equations we have
\begin{align}
 & \left| \int_{y > Y_0} y f''(x) \chi (y) \Im[ s_N (x+ \i y ) - s(x+ \i y ) ] \d x \d y \right| \notag\\
 = & \left| \int_{y > Y_0} y f'(x) \chi (y) \Im[ \del_x (s_N (x+ \i y ) - s(x+ \i y ) )] \d x\d y \right| \notag\\
 \leq & \left| \int_{2A> y > Y_0} |y f'(x)| | \del_z (s_N (x+ \i y ) - s(x+ \i y ) ) | \d x \d y \right|
\end{align}
By Cauchy's integral formula we can write
\beq
\del_z (s_N (z) - s(z) ) = \frac{1}{2 \pi \i} \int_{ w : |w-z| = \Im[z]/2} \frac{ s_N(w) -s (w) }{(w-z)^2} \d w .
\eeq
Therefore by Minkowski's integral inequality and \eqref{eqn:BMP-s} we have
\beq
\| \del_z (s_N (z) - s(z)) \|_p \lesssim \frac{1}{ (\Im[z] )^2} .
\eeq
Therefore,
\begin{align}
\left\| \int_{2A> y > Y_0} |y f'(x)| | \del_z (s_N (x+ \i y ) - s(x+ \i y ) ) | \d x \d y \right\|_p \lesssim \| f'\|_{L^1} \int_{Y_0 < y< 2A} \frac{y}{y^2} \d y \lesssim \log (E_0) .
\end{align}
This completes the proof. \qed

We define the eigenvalue counting functions by
\beq
\N (E) := | \{ i : x_i \leq E \} |, \qquad \nsc (E) := \int_{-\infty}^E \rho (x) \d x.
\eeq
\bel \label{lem:rig-N}
We have for $10 \leq E \leq c_* N^{2/3}$,
\beq
\| \N(E) - \nsc (E) \|_p \lesssim \log (E)
\eeq
\eel
\proof Let $f_1$ be a smooth function s.t. $f_1(x) =1$ for $x \in [-100, E]$ and $f_1(x) = 0$ for $x > E + (E)^{-1/2}$ or $x < -100 - E^{-1/2}$. Let $f_2 (x)$ be a function s.t. $f_2(x) = 1$ for $x\in [-100, E-E^{-1/2} ]$ and $f_2 (x) = 0$ for $x < -100 - E^{-1/2}$ or $x > E$. We further assume $0 \leq f_i \leq 1$. We have
\beq
|\nsc (E) - \int f_i (x) \rho (x) \d x | \lesssim 1
\eeq
for $i=1, 2$. It is therefore easy to see that
\beq
| \N (E)  - \nsc(E) | \lesssim \sum_{i=1}^2 \left| \sum_j f_i (x_j) - \int f_i (x) \rho (x) \d x \right| + 1 + \N (0) .
\eeq
From \cite[Lemma 3.6]{bourgade2022optimal} we have
\beq
\| \N (0) \|_p \lesssim 1. 
\eeq
The claim follows from the above two inequalities and a direct application of Lemma \ref{lem:rig-proof-1}. \qed

\vspace{5 pt}

\noindent{\bf Proof of Proposition \ref{prop:rigidity}}.  The estimate \eqref{eqn:usual-rigid} is well-known and follows from \cite[Theorem 2.4]{bourgade2014edge}. We turn to the proof of \eqref{eqn:rigid-moment}. Throughout we assume $j \leq N^{1/10}$.  By definition
\beq
\int_{0}^{\tilgam_j} \rho (x) \d x = j.
\eeq
It is not hard to see that for $0 \leq s \leq \tilgam_j$ and $j \leq N/2$ we have,
\beq \label{eqn:rig-proof-a2}
\nsc (\tilgam_j) - \nsc (s) \asymp (\tilgam_j - s) j^{1/3} 
\eeq
and for $\tilgam_j \leq s \leq 4 N^{2/3}$ that
\beq \label{eqn:rig-proof-a3}
\nsc (s) - \nsc (\tilgam_j) \asymp (s- \tilgam_j ) s^{1/2} .
\eeq
Let $K \geq 10$ be the smallest $K$ s.t. $\tilgam_K \geq 100$. For $0 \leq s \leq \tilgam_j - 10$ we have for $j \geq K$ that
\begin{align}
\pp\left[ x_j < \tilgam_j - s \right] &= \pp\left[ \N( \tilgam_j -s ) - \nsc ( \tilgam_j - s) \geq \nsc (\tilgam_j) - \nsc (\tilgam_j - s) \right]  \notag\\
\leq & C_p \left( \frac{ \log j}{ s j^{1/3} } \right)^p 
\end{align}
by Markov, Lemma \ref{lem:rig-N} and \eqref{eqn:rig-proof-a2}. By adjusting constants one can improve this to all $0 \leq s \leq 2 \tilgam_j$, and then by applying \cite[Corollary 1.6]{bourgade2022optimal} to all $0 \leq s \leq c_* N^{2/3}$.  Since these estimates hold for any $p \geq 1$, one finds that
\beq
\ee\left[ |x_j - \tilgam_j |^p \1_{ \{- c_* N^{2/3} \leq x_j \leq \tilgam_j\} } \right] \lesssim \left( \frac{ \log j }{j^{1/3}} \right)^p 
\eeq
It is not hard to derive the estimate
\beq
\ee[ |x_j - \tilgam_j|^p \1_{ \{ x_j \leq - c_* N^{2/3} \} } ] \lesssim N^{-D}
\eeq
for any $D \geq 1$ 
using \cite[Corollary 1.6]{bourgade2022optimal} and a weak estimate $\ee[ |x_1|^{2p} ] \leq N^{20 p}$ (which can be derived, e.g., from the tri-diagonal representation of the Gaussian $\beta$ ensemble). In a similar manner we obtain
\beq
\pp\left[ x_j > \tilgam_j + s \right] \leq C_p \left( \frac{\log j}{ s j^{1/3} } \right)^p
\eeq
for all $s \leq c_* N^{2/3}$. Proceeding as above we therefore obtain
\beq
\ee\left[ |x_j - \tilgam_j |^p \right] \leq C_p \left( \frac{ \log j}{j^{1/3}} \right)^p
\eeq
for all $j \geq K$. For $j \leq K$,  we can bound $x_j -\tilgam_j \leq |x_K - \tilgam_K| + \O_K (1)$ as well as $x_j - \tilgam_j \geq x_1 - \tilgam_1 - \O_K (1)$. The moments of the first quantity we have already bounded. The moments of the latter quantity can be bounded using again \cite[Corollary 1.6]{bourgade2022optimal}. This completes the proof of \eqref{eqn:rigid-moment} after noting that $| \hatgam_i - \tilgam_i| \leq C i^{-1/3}$ for $i \leq N^{2/5}$ using \eqref{eqn:quantiles-approx}. \qed

\bel \label{lem:aux-rigidity} Let $x_i$ be a rigid process. We have with overwhelming probability,
\beq
\sup_{ |y| \leq 2N^{1/100} } \left| \sum_{i= N^{1/10}}^N \frac{ z_i }{ y -x_i} + \int_{\tilgam_{N^{1/10}}}^\infty \frac{ \rho (x) }{ x } \d x \right| \leq N^{-1/100}
\eeq
\eel
\proof Since $x_{N^{1/10} } \geq N^{1/20}$ with overwhelming probability, we see that the LHS is a Lipschitz function with constant bounded above by $N$ with overwhelming probability on the domain $[-2N^{1/100} , 2N^{1/100} ]$ so it suffices to prove the estimate for fixed $y$. It is easy to see that with overwhelming probability, using $|z_i| \leq N^{1/1000}$ with  overwhelming probability, and rigidity for $x_i$
\beq
\left| \sum_{i=N^{1/10}}^N \frac{z_i}{ y-x_i} + z_i \int_{\tilgam_{i-1}}^{\tilgam_i} \frac{\rho (x) }{x} \right| \lesssim N^{1/1000} \sum_{i=N^{1/10}}^N \frac{ N^{1/100} }{i^{4/3}} \leq N^{ -1/90} .
\eeq
On the other hand by  standard large deviations  bounds (see, e.g., \cite[Appendix C]{erdHos2013local}) it is easy to see that
\beq
\left| \sum_{i=N^{1/10}}^N (z_i -1 ) \int_{\tilgam_{i-1}}^{\tilgam_i} \frac{ \rho (x) }{x} \d x \right| \leq N^{-1/90}
\eeq
with overwhelming probability. This yields the claim. \qed

\subsection{Proof of Lemma \ref{lem:iso}}
\label{sec:iso}
By standard large deviations bounds \cite[Appendix C]{erdHos2013local} and rigidity \eqref{eqn:usual-rigid} we have that with overwhelming probability,
\beq
\left| \sum_{i=1}^N \frac{z_i -1}{  E - x_i} \right| \leq N^{\eps} \left( \sum_{i=1}^N \frac{1}{ i^{4/3} + E^2} \right)^{1/2} \lesssim N^{\eps} |E|^{-1/4} ,
\eeq
for any $\eps >0$. With overwhelming probability,
\begin{align}
 & \left| \sum_{i=1}^N \frac{1}{ E - x_i}  - \int_{0}^\infty \frac{ \rho (x) }{E-x} \d x \right| \leq \sum_{i=1}^N \int_{\tilgam_{i-1}}^{\tilgam_i} \left| \frac{1}{ x_i - E} - \frac{1}{ x - E} \right| \rho (x) \d x \notag\\
 \lesssim & N^{\eps} \sum_{i=1}^N \frac{1}{\min \{ i^{1/3} , ( N +1 - i )^{1/3} \}} \frac{1}{i^{4/3} + E^2} \lesssim \frac{N^{\eps}}{|E|^{1/2} } . 
\end{align}
This yields the claim. \qed

\section{Deterministic lemma}

\bel \label{lem:calc-1}
Let $M \geq 1$ and $m \geq 1$. There is a constant $c > 0$ so that for all $M \leq m^{2/3}$ we have
\beq
\int_0^{\hatgam_m} \sqrt{x} \left( \frac{1}{x} - \frac{1}{x+M} \right) \d x \geq  c M^{1/2} 
\eeq
\eel
\proof The integrand is positive so, for some $c' > 0$, we may bound it below by
\begin{align}
\int_0^{c' M} \sqrt{x} \left( \frac{1}{x} - \frac{1}{x+M} \right) \d x = \int_0^{c' M} \frac{M \sqrt{x}}{x(x+M)} \d x \geq \frac{1}{c' + 1} \int_0^{c'M} \frac{1}{\sqrt{x} } \d x \geq c M^{1/2}
\end{align}
for some $c > 0$, which yields the claim. \qed



\bibliography{mybib}{}

\begin{thebibliography}{10}

\bibitem{adler2013random}
M.~Adler, P.~van Moerbeke, and D.~Wang.
\newblock Random matrix minor processes related to percolation theory.
\newblock {\em Random Matrices: Theory and Applications}, 2(04):1350008, 2013.

\bibitem{baik2005phase}
J.~Baik, G.~Ben~Arous, and S.~P{\'e}ch{\'e}.
\newblock Phase transition of the largest eigenvalue for nonnull complex sample
  covariance matrices.
\newblock {\em Annals of Probability}, 33(5):1643--1697, 2005.

\bibitem{bao2022eigenvector}
Z.~Bao and D.~Wang.
\newblock Eigenvector distribution in the critical regime of {BBP} transition.
\newblock {\em Probability Theory and Related Fields}, 182(1):399--479, 2022.

\bibitem{bao2026eigenvectordistributionrandommatrices}
Z.~Bao, D.~Wang, and Y.~Zhu.
\newblock Eigenvector distribution of random matrices under critical
  finite-rank deformations.
\newblock {\em Preprint, arXiv:2605.30779}, 2026.

\bibitem{alex2013limits}
A.~Bloemendal and B.~Vir{\'a}g.
\newblock Limits of spiked random matrices {I}.
\newblock {\em Probability Theory and Related Fields}, 156(3-4):795--825, 2013.

\bibitem{bloemendal2016limits}
A.~Bloemendal and B.~Vir{\'a}g.
\newblock Limits of spiked random matrices {II}.
\newblock {\em Annals of Probability}, 44(4):2726--2769, 2016.

\bibitem{bourgade2014edge}
P.~Bourgade, L.~Erd{\"o}s, and H.-T. Yau.
\newblock Edge universality of beta ensembles.
\newblock {\em Communications in Mathematical Physics}, 332(1):261--353, 2014.

\bibitem{bourgade2022optimal}
P.~Bourgade, K.~Mody, and M.~Pain.
\newblock {Optimal local law and central limit theorem for $\beta$-ensembles}.
\newblock {\em Communications in Mathematical Physics}, 390(3):1017--1079,
  2022.

\bibitem{bykhovskaya2026weakweakfactorsuniform}
A.~Bykhovskaya, V.~Gorin, and S.~Sodin.
\newblock How weak are weak factors? {U}niform inference for signal strength in
  signal plus noise models.
\newblock {\em Preprint, arXiv:2507.18554}, 2026.

\bibitem{erdHos2013local}
L.~Erd{\H{o}}s, A.~Knowles, H.-T. Yau, and J.~Yin.
\newblock The local semicircle law for a general class of random matrices.
\newblock {\em Electronic Journal of Probability}, 18:1--58, 2013.

\bibitem{erdHos2017dynamical}
L.~Erd{\H{o}}s and H.-T. Yau.
\newblock {\em A dynamical approach to random matrix theory}, volume~28.
\newblock American Mathematical Society, 2017.

\bibitem{forrester2010log}
P.~J. Forrester.
\newblock {\em Log-gases and random matrices (LMS-34)}.
\newblock Princeton University Press, 2010.

\bibitem{gorin2015multilevel}
V.~Gorin and M.~Shkolnikov.
\newblock {Multilevel Dyson Brownian motions via {J}ack polynomials}.
\newblock {\em Probability Theory and Related Fields}, 163(3):413--463, 2015.

\bibitem{gorin2014interactingparticlesystemsedge}
V.~Gorin and M.~Shkolnikov.
\newblock {Interacting particle systems at the edge of multilevel Dyson
  Brownian motions}.
\newblock {\em Advances in Mathematics}, 304:90--130, 2017.

\bibitem{gustavsson2005gaussian}
J.~Gustavsson.
\newblock {Gaussian fluctuations of eigenvalues in the GUE}.
\newblock {\em Annales de l'Institut Henri Poincar{\'e} (B) Probabilit{\'e}s et
  statistiques}, 41(2):151--178, 2005.

\bibitem{huang2019eigenvalues}
J.~Huang.
\newblock Eigenvalues for the minors of {W}igner matrices.
\newblock {\em Annales de l'Institut Henri Poincar{\'e} (B) Probabilit{\'e}s et
  statistiques}, 58(4):2201--2215, 2022.

\bibitem{landon2022fluctuations}
B.~Landon and P.~Sosoe.
\newblock {Fluctuations of the overlap at low temperature in the 2-spin
  spherical SK model}.
\newblock {\em Annales de l'Institut Henri Poincar{\'e} (B) Probabilit{\'e}s et
  statistiques}, 58(3):1426--1459, 2022.

\bibitem{landon2020fluctuations}
B.~Landon and P.~Sosoe.
\newblock {Fluctuations of the 2-spin SSK model with magnetic field}.
\newblock {\em Annals of Applied Probability}, 35(5):3543--3619, 2025.

\bibitem{najnudel2021bead}
J.~Najnudel and B.~Vir{\'a}g.
\newblock The bead process for beta ensembles.
\newblock {\em Probability Theory and Related Fields}, 179(3):589--647, 2021.

\bibitem{peche2006largest}
S.~P{\'e}ch{\'e}.
\newblock The largest eigenvalue of small rank perturbations of hermitian
  random matrices.
\newblock {\em Probability Theory and Related Fields}, 134(1):127--173, 2006.

\bibitem{petrov2026perturbedbetacornersprocess}
L.~Petrov and J.~Xu.
\newblock Perturbed beta corners process.
\newblock {\em Preprint, arXiv:2607.27810}, 2026.

\bibitem{ramirez2011beta}
J.~Ram{\'{\dotlessi}}rez, B.~Rider, and B.~Vir{\'a}g.
\newblock {Beta ensembles, stochastic Airy spectrum, and a diffusion}.
\newblock {\em Journal of the American Mathematical Society}, 24(4):919--944,
  2011.

\bibitem{tracy1994level}
C.~A. Tracy and H.~Widom.
\newblock Level-spacing distributions and the {A}iry kernel.
\newblock {\em Communications in Mathematical Physics}, 159(1):151--174, 1994.

\bibitem{tracy1996orthogonal}
C.~A. Tracy and H.~Widom.
\newblock On orthogonal and symplectic matrix ensembles.
\newblock {\em Communications in Mathematical Physics}, 177(3):727--754, 1996.

\end{thebibliography}
\bibliographystyle{abbrv}

\end{document}